\documentclass[10pt]{article}
\usepackage{enumerate, fullpage}
\usepackage{graphicx}
\usepackage{multirow}
\usepackage{color}
\usepackage{amsmath}
\usepackage{amsfonts}
\usepackage{amssymb}
\usepackage{mathrsfs}
\usepackage{eucal,bm}
\usepackage{lscape}
\usepackage{mathptmx}
\definecolor{darkgreen}{RGB}{0, 102, 0}
\definecolor{IKB}{RGB}{0, 47, 167}
\usepackage[pdftex,bookmarks,bookmarksopen=true,pdfstartview={FitH},bookmarksnumbered=true,colorlinks=true,citecolor=darkgreen,linkcolor=IKB,plainpages=false,pdfpagelabels]{hyperref}

\def\norm#1{\|#1\|}

\newcommand{\R}{\mathbb{R}}

\newcommand{\beq}{\begin{equation}}
\newcommand{\eeq}{\end{equation}}

\newcommand{\beqnr}{\begin{eqnarray}}
\newcommand{\eeqnr}{\end{eqnarray}}

\newcommand{\benum}{\begin{enumerate}}
\newcommand{\eenum}{\end{enumerate}}
\newcommand{\argmax}{\mathop{\rm argmax}}

\newcommand{\sgn}{\mathop{\rm sgn}}
\newcommand{\QED}{\rule{.1in}{.1in}}

\newtheorem{DE}{Definition}[section]
\newtheorem{AS}[DE]{Assumption}

\newtheorem{LE}[DE]{Lemma}

\newtheorem{RE}[DE]{Remark}
\newtheorem{THM}[DE]{Theorem}
\newtheorem{CO}[DE]{Corollary}

\newcommand{\qed}{\mbox{}\hspace*{\fill}\nolinebreak\mbox{$\rule{0.7em}{0.7em}$}}

\advance \topmargin by -\headheight
\advance \topmargin by -\headsep

\evensidemargin \oddsidemargin

\begin{document}

\begin{center}
\Large {\bf A note on polynomial solvability of the CDT problem}

\vskip 0.1cm

\normalsize

Daniel Bienstock, Columbia University, December 2013
\end{center}

\begin{abstract}
Recent research has focused on the complexity of extensions of the 
classical trust-region subproblem, which addresses the minimization of a 
quadratic function over a unit ball in $\R^n,$ a problem of importance in 
many applications.   Even though the 
trust-region subproblem can be considered well-solved (both from the perspective of theory and implementation) even minor extensions are NP-hard.
The CDT (Celis-Dennis-Tapia) problem is an extension of the trust-region subproblem, involving the minimization of a quadratic  
function over the intersection of two ellipsoids in $\R^n$.  In this paper
we show how to adapt a construction of Barvinok so as to obtain a polynomial-time algorithm for quadratic programming with
a fixed number of quadratic constraints (one of which is ellipsoidal) under the bit model of computing.
\end{abstract}

\section{Introduction}\label{intro}

The classical trust-region subproblem can be stated as $\min \{ f(x) \, : \, \norm{x} \le 1, \, x \in \R^n \}$ where $f \, : \, \R^n \rightarrow \R$ is a 
quadratic function, i.e. a polynomial of degree at most two.  This problem, which has received a great deal of attention
in the literature because of a rich set of applications, can be efficiently solved -- it is polynomial-time solvable and effective implementations are
available.  See \cite{trustbook} (and references therein), \cite{yeold}.  At
the same time, seemingly minor extensions of the trust-region subproblem, such
as the addition of an arbitrary family of linear constraints, give rise to
NP-hard problems, and recent work has focused on understanding which
extensions are polynomial-time solvable.  See \cite{bureranstreicher}, \cite{yangburer0}, \cite{danoalex} for recent complexity results.

The CDT (Celis-Dennis-Tapia) problem \cite{CDT} is of the form
$$(CDT) \quad \min \{ f(x) \, : \, g_2(c) \le 0, \ g_2(x) \le 0, \, x \in \R^n \}$$ where $f(x)$ is a quadratic and for $i = 1, 2$ the set $\{x \in \R^n \, : \, g_i(x) \le 0 \}$ defines an
ellipsoid. Thus, the CDT problem can be viewed as an extension of the trust-region subproblem with an added ellipsoidal constraint.  The CDT problem has long generated
interest, see \cite{bureranstreicher}, where it was termed the ``two trust-region'' problem, \cite{aizhang}, \cite{beckeldar}, \cite{cy99}, \cite{cy09}, \cite{pengyuan}, \cite{yz}, \cite{JeLeeLi}, \cite{pengyuan}, \cite{yuan}; also see \cite{yangburer} and \cite{bomzeoverton}, as well as their references.  Broadly speaking these papers have sought to exploit 
the connection between CDT and convex optimization, in particular semidefinite programming; 
this approach is related to the use of the S-Lemma to solve
the classical trust-region subproblem (see \cite{polter}).
The inclusion of a second ellipsoidal constraint dramatically increases the complexity of the underlying geometry, even in low dimension, as demonstrated in particular in the elegant
analysis in \cite{yangburer}. From our perspective, this analysis highlights an underlying combinatorial (or perhaps, enumerational) behavior of
CDT problems as evinced by the local structure of solutions to systems of quadratic inequalities. Recently,  an 
algorithm for quadratic programming with two quadratic constraints, which runs in polynomial time under appropriate assumptions, was
presented in \cite{SNTI}.

A separate line of work has produced results of a very different flavor in order to 
address related problems.  Barvinok \cite{barvinok}
(``Problem(1.1)'')
proved that for each each positive integer $K$ there is a polynomial-time  algorithm to decide if a system of the form
\begin{subequations}
\label{barv}
\begin{eqnarray}
x^T M_i x & = & 0, \quad 1 \le i \le K, \\
\| x \|^2 & = & 1,
\end{eqnarray}
 \end{subequations}
is feasible, where $x \in \R^n$ and $M_i$ is an $n \times n$ matrix for $1 \le i \le K$.  The construction in \cite{barvinok} relies on the model of computing over the reals (e.g., infinite-precision is assumed) however the extension 
to the standard bit model of computing should be straightforward \cite{barvpersonal}.   

Related and stronger results were presented by Grigoriev and Pasechnik in \cite{dima}.  Moreover, it was argued in \cite{dima2} that the results
in \cite{dima} imply that a polynomial-time algorithm for a generalization
of CDT exists.  Also see \cite{baspasu}.
Yet another line of research comes from the algebraic and computational geometry
communities, in particular the theory of ``roadmaps'' of semi-algebraic sets, which 
was started in \cite{canny}.   This topic appears  related to the work cited above on ``sampling'' algebraic sets. See \cite{basuetal}, and \cite{seldinschost} for
recent research results and additional citations.  These research efforts have 
produced results with applications to diverse
problem domains. It is quite possible that other polynomial-time algorithms for generalizations of CDT 
can be derived from this work.  A nontrivial point concerning such
algorithms is whether they would attain
polynomial-time complexity under the bit model of computation, as opposed to exact
computation over the reals.

In this note we present a simple procedure that uses a weak version of Barvinok's construction
to obtain a polynomial-time algorithm for a generalization of CDT under the standard bit model of computing.   This algorithm will use a relaxed version of feasibility, as follows.
\begin{DE}\label{epsfeas}  Let $f_i(x) \le 0$, $1 \le i \le m$ be a system
of inequalities where $x \in \R^n$.
\begin{itemize}
\item [(a)]  Given $0 \le \epsilon < 1$ a vector $\hat x \in \R^n$ is called $\epsilon$-feasible for the system if $f_i(\hat x) \le \epsilon$ for $1 \le i \le m$.  If such a vector exists we will say that the system
is $\epsilon$-feasible.  
\item [(b)] An algorithm that, given any $0 < \epsilon < 1$ either proves that the system is
infeasible, or proves that it is $\epsilon$-feasible will be 
called a  {\bf weak feasibility algorithm}.  Note that we insist on $\epsilon > 0$.
\end{itemize}
\end{DE}
The weak version of Barvinok's construction is as follows:
\begin{AS} \label{weakalgo}We will assume that for each integer $K$, there is a
weak feasibility algorithm for systems of the type (\ref{barv}) 
with running time polynomial in the size of the data and in $\log \epsilon^{-1}$.
\end{AS}
Barvinok's method \cite{barvinok} clearly fulfills this role, as does Grigoriev and Pasechnik's
\cite{dima}. It is also possible that faster weak feasibility algorithms exist, as opposed to algorithms for exact
feasibility as in (\ref{barv}).  Also see \cite{barvlinf}.\\

To describe our main result suppose that for $0 \le i \le p$, $\ g_i(x)$ is a quadratic, over $x \in \R^n$.   We consider the problem
\begin{eqnarray}
&&  \min \, \{ \, g_0(x) \, : \, g_i(x) \le 0, \ 1 \le i \le p \, \}, \label{cdtgen}
\end{eqnarray}
and prove:
\begin{THM} \label{main} For each fixed integer $p \ge 1$ there is an algorithm with the following properties. Given a problem of the form (\ref{cdtgen})
where at least one of the $g_i(x)$ with $i \ge 1$ is strictly convex, and $0 < \epsilon < 1$, 
the algorithm either 
\begin{itemize}
\item [(1)] proves
that problem (\ref{cdtgen}) is infeasible, 
\end{itemize}
or
\begin{itemize}
\item [(2)] computes an $\epsilon$-feasible vector $\hat x$ such that there exists no feasible $x$ with $g_0(x) < g(\hat x) - \epsilon$. 
\end{itemize}
Under Assumption \ref{weakalgo} the algorithm runs in polynomial time. 
More precisely, the algorithm makes a sequence of calls to a weak feasibility algorithm 
for problems of type (\ref{barv}) with $K = O(p)$; 
the length of the sequence is polynomial in the number of bits in the
data and $\log \epsilon^{-1}$, as is the size of the coefficients 
of the matrices $M_i$, and as is all additional work carried out by the algorithm.
\end{THM}
Thus, under Assumption \ref{weakalgo}, Theorem \ref{main} implies that a polynomial-time algorithm
for CDT exists. Theorem \ref{main} is proved in several steps in Section \ref{construction}.  
 Throughout, we assume $n \ge 2$.
\section{The construction}\label{construction}
This section is organized as follows.
In Section \ref{algorithm1} we describe an algorithm to determine if a 
system of $m$ quadratic inequalities is $\epsilon$-feasible; this algorithm runs
in polynomial time  for each fixed $m$ provided that at least one of the 
quadratics is strictly convex.
In Section \ref{binsearch} 
the algorithm in Section \ref{algorithm1} is used to compute the value of problem (\ref{cdtgen}) 
within tolerance $\epsilon$, in polynomial time, under the assumptions in 
Theorem \ref{main}.
However, this does not yet yield a proof of Theorem \ref{main} because the  algorithm we describe in Section \ref{algorithm1} relies on the weak feasibility
algorithm in Assumption \ref{weakalgo} as a subroutine.
That algorithm (and in the strict sense, Barvinok's) decides if a system of the form (\ref{barv}) is $\epsilon$-feasible
but without producing an explicit $\epsilon$-feasible vector. 
In Section \ref{explicit} we show how to refine  our algorithm
from Section \ref{algorithm1} so as to produce an $\epsilon$-feasible vector in the case that infeasibility is not proved.  This fact, together with the results
in Section \ref{binsearch} will be used to complete the proof of Theorem \ref{main}.

\subsection{Systems of quadratic inequalities}\label{algorithm1}
Here we consider a system of quadratic inequalities
\begin{eqnarray}
f_i(x) & \le & 0, \quad 1 \le i \le m. \label{manyq}
\end{eqnarray}
We write
\begin{eqnarray}
f_i(x) & \doteq & x^T A_i x + c^T_i x + d_i \label{quad-def} 
\end{eqnarray}
where $A_i \in \R^{n \times n}$ and symmetric, $c_i \in \R^n$ and $d_i \in \R$. 
Such a system describes
the feasibility set for a problem of the form (\ref{cdtgen}); more generally
we will use the solution of systems of the form (\ref{manyq}) as steps in
our algorithm for problem (\ref{cdtgen}).  Our main result for this section, proved below, is
as follows:

\begin{THM} \label{punch} Under Assumption \ref{weakalgo}, for each fixed $m$ there is a polynomial-time 
weak feasibility algorithm for  any system of type (\ref{manyq}) 
where $A_i \succ 0$ for at least one index $i \ge 1$. \end{THM}

\noindent We will first prove two technical results, Lemma \ref{equiv0} and
\ref{equiv}, under the assumptions of Theorem \ref{punch}. We assume, without loss of generality that
\begin{eqnarray}
 f_1(x) & = & \norm{x}^2 - 1. \label{one}
\end{eqnarray}
Then, for $2 \le i \le m$, there exists (polynomially computable) $U_i > 0$, 
such that
$$ |f_i(x)| \le U_i, \quad \mbox{for each $x \in \R^n$ with $\norm{x}^2 \le 2$.}$$
Now consider the following system of quadratic equations on 
real variables $v_0, x_1, \ldots, x_n, s_1, \ldots, s_m, w_2, \ldots, w_m$:
\begin{subequations}
\label{sys2}
\begin{eqnarray}
x^T A_i x + c^T_i v_0 x + d_i v_0^2 + s_i^2 & = & 0 \quad \quad 1 \le i \le m, \label{sys2a} \\
\frac{s_i^2 + w_i^2}{U_i} - v_0^2 & = & 0 \quad \quad 2 \le i \le m,\label{sys2b} \\
\norm{x}^2 + s_1^2 + \sum_{i = 2}^n \frac{s_i^2 + w_i^2}{U_i} + v_0^2 
& = & m + 1. \label{sys2c}
\end{eqnarray}
\end{subequations}

\begin{LE} \label{equiv0}  Let $0 \le \delta < 1$ and suppose that
$$\hat z \ = \ (\hat v_0, \hat x_1, \ldots, \hat x_n, \hat s_1, \ldots, \hat s_m, \hat w_2, \ldots, \hat w_m)^T$$
is a 
$\delta$-feasible solution to (\ref{sys2a})-(\ref{sys2b}).  Then  {\bf(i)}
\begin{eqnarray}
&& m \hat v_0^2 - m \delta \ \le \ \norm{\hat x}^2 + \hat s_1^2 + \sum_{i = 2}^n \frac{\hat s_i^2 + \hat w_i^2}{U_i} \ \le \ m \hat v_0^2 + m \delta. \label{thesum}
\end{eqnarray}
{\bf(ii)} If $\hat z$ is also $\delta$-feasible for (\ref{sys2c})
$$1 - \delta \le \hat v_0^2 \le  1 + \delta.$$ \end{LE}
\noindent {\em Proof.} (i) Note that $\delta$-feasibility of $\hat x$ applied to (\ref{sys2a}) for $i = 1$ 
states
$$ -\delta \ \le \ \norm{\hat x}^2 - \hat v_0^2 + \hat s^2_1 \ \le \ \delta, $$
and applied to (\ref{sys2b}) it states
$$ -\delta \ \le \ \frac{\hat s_i^2 + \hat w_i^2}{U_i} - \hat v_0^2 \ \le \ \delta, \quad 2 \le i \le m.$$
Adding these inequalities yields (\ref{thesum}).
(ii) Since $\hat x$ is $\delta$-feasible for (\ref{sys2c}) 
$$ m+1 - \delta \ \le \ \norm{\hat x}^2 + \hat s_1^2 + \sum_{i = 2}^n \frac{\hat s_i^2 + \hat w_i^2}{U_i} + \hat v_0^2 \ \le \ m+1 + \delta$$
which together with (\ref{thesum}) implies the desired result. \QED \\

\noindent Let $M$ denote the largest absolute value of a coefficient in system (\ref{sys2}).

\begin{LE} \label{equiv}  Let $ 0 \le \epsilon < 1/2$.  (a) Suppose (\ref{manyq}) is $\epsilon$-feasible.  Then
(\ref{sys2}) is $m\epsilon$-feasible.   (b) Conversely, 
if (\ref{sys2}) is $\epsilon$-feasible, then (\ref{manyq}) is $(2n + 1)M \epsilon$-feasible. \end{LE}
\noindent {\em Proof.} (a) Suppose first that $\tilde x$ is an $\epsilon$-feasible solution to (\ref{manyq}). 
Define:
\begin{eqnarray}
\tilde v_0 & \doteq & 1,\\
\tilde s_i & \doteq & \sqrt{ \max\{0, -f_i(\tilde x)\}}, \quad 1 \le i \le m, \label{slack}\\
\tilde w_i & \doteq & \sqrt{ U_i - \tilde s_i^2}, \quad 2 \le i \le m. \label{comp}
\end{eqnarray}
In (\ref{comp}), the value in the radical is nonnegative  since $|\tilde x|^2 \le 1 + \epsilon < 2$. Now
we claim that 
$$ \tilde z \doteq (\tilde x_1 , \ldots, \tilde x_n, \tilde v_0, \tilde s_1, \ldots, \tilde s_m, \tilde w_2, \ldots, \tilde w_m)^T$$ 
is an $\epsilon$-feasible solution to (\ref{sys2}). To see this, note that
$\tilde v_0 = 1$ and (\ref{slack}) imply that $\tilde z$ is $\epsilon$-feasible
for (\ref{sys2a}). Likewise, $\tilde z$ satisfies (\ref{sys2b}) by (\ref{comp}) and $\tilde v_0 = 1$.  Finally, Lemma \ref{equiv0}, part (i) (with $\hat x = \tilde x$ and
$\delta = \epsilon$), together with $\tilde v_0 = 1$ implies that
$\tilde z$ is $m\epsilon$-feasible for (\ref{sys2c}), as desired. \\

\noindent (b) For the converse, suppose 
$$ \hat z \doteq (\hat x_1 , \ldots, \hat x_n, \hat v_0, \hat s_1, \ldots, \hat s_m, \hat w_2, \ldots, \hat w_m)^T$$
is $\epsilon$-feasible 
for  (\ref{sys2}). By Lemma \ref{equiv0} (ii), 
\begin{eqnarray}
&& 1 - \epsilon \ \le \ \hat v_0^2 \le 1 + \epsilon. \label{boundonvprev}
\end{eqnarray}
This implies
\begin{eqnarray}
&& 1 - \epsilon \ \le \ | \hat v_0| \le 1 + \epsilon, \label{boundonv0}
\end{eqnarray}
and together with is $\epsilon$-feasibility of $\hat x$ 
for (\ref{sys2a}) with $i = 1$, it implies $ \norm{\hat x}^2 + \hat s_1 ^2 \ \le \ 1 + 2 \epsilon $
and therefore
\begin{eqnarray}
&& | \hat x_j | \le 1 + \epsilon, \quad \mbox{for $1 \le j \le n$}.  \label{boundonxj}
\end{eqnarray}
Then (\ref{boundonvprev})-(\ref{boundonxj}) imply that for each $1 \le j \le n$, $|v_0 \hat x_j - \hat x_j| \le \epsilon |\hat x_j| \le \epsilon^2 + 2 \epsilon < 2$. Hence, using ``$\sgn$'' indicate the sign function,
 $\epsilon$-feasibility of $\hat x$ for (\ref{sys2a}) imply that 
$\sgn(\hat v_0) (\hat x_1 , \ldots, \hat x_n)^T$ is  $(2n + 1)M \epsilon$-feasible
for (\ref{manyq}). \QED

\begin{CO} \label{approx} A system
(\ref{manyq}) is feasible if and only if the corresponding system (\ref{sys2}) is feasible. \end{CO}
\noindent {\em Proof.} Use $\epsilon = 0$ in Lemma \ref{equiv}.  \QED \\

\noindent Now we can present the proof of the main result in this section.\\

\noindent {\em Proof of Theorem \ref{punch}.}  Consider the corresponding system (\ref{sys2}),
and let $\delta \doteq \frac{\epsilon}{(2n + 1)M}$. 
Using the method in Assumption \ref{weakalgo}, we terminate in polynomial time
with a proof that system  (\ref{sys2}) is infeasible, in 
which case system (\ref{manyq}) is infeasible (by Corollary \ref{approx}), 
or, using part (b) of Lemma \ref{equiv}, with a proof that (\ref{manyq}) is
$\epsilon$-feasible. \QED

\subsection{Estimating the value of problem (\ref{cdtgen})}\label{binsearch}
In what follows we will use the following convention.
\begin{DE} Consider an optimization problem 
\begin{eqnarray}
&& H^* \ \doteq \ \min \{ f(x) \, : \, h_i(x) \le 0, \ 1 \le i \le k \}. \label{Hopt}
\end{eqnarray}
Given $0 < \epsilon < 1$ a rational $V$ is called an $\epsilon$-estimate for (\ref{Hopt}) if $V \le H^*$ and there exists some vector $\bar x$ with $h_i(x) \le \epsilon$ for 
$1 \le i \le k$ and $f(\bar x) \le V + \epsilon$.  We will likewise
define $\epsilon$-estimates for maximization problems.
\end{DE}
Let $G^*$ be the value of problem (\ref{cdtgen}), i.e. $G^* \ \doteq \ \min \{ g_0(x) \, : \, g_i(x) \le 0, \ 1 \le i \le p \}.$
In this section we provide a polynomial-time procedure that produces one of
two outcomes, given $0 < \epsilon < 1:$
\begin{itemize}
\item[(i)] It proves that (\ref{cdtgen}) is infeasible.
\item [(ii)] It produces an $\epsilon$-estimate  for (\ref{cdtgen}).
\end{itemize}
We stress that in case (ii) the vector $\bar x$ that yields the $\epsilon$-estimate is not explicitly known. The algorithm, which amounts to a modified form of binary search, works as follows. As the initial step we run the weak feasibility problem on the 
system $\{ x \in \R^n : \, g_i(x) \le 0, \ 1 \le i \le p \},$ and if infeasibility is determined, we stop with outcome (i).  Otherwise we know some
$\epsilon$-feasible point exists. We next compute a rational $U > 1$ such that
$ |g_0(x)| \le U $ for every $\epsilon$-feasible $x$.  Such a value $U$ 
exists
and is polynomial-time computable since we assume that at least 
one of the $g_i(x)$, for $i \ge 1$, is positive definite. 

In a typical iteration of the modified binary search we solve the 
weak feasibility problem, with tolerance $\epsilon$, for the system of the form
\begin{subequations}
\label{redx}
\begin{eqnarray}
g_i(x) & \le & \ \, 0  \quad \mbox{for} \ 1 \le i \le p, \\
g_0(x) & \le &  \ \, V 
\end{eqnarray}
\end{subequations}
If (\ref{redx}) is infeasible, then clearly $V \le G^*$ and otherwise 
there exists an $\epsilon$-feasible $\bar x$ with $g_0(\bar x) \le V + \epsilon$.  Clearly, 
after a number of iterations which is $O(\log U + \log \epsilon^{-1})$ (i.e. is polynomial in the number of bits in the data) we will terminate with one of
the two desired outcomes.\\

\noindent {\bf Remark.} Assume (\ref{cdtgen}) is feasible. For a given $\epsilon$ denote by $\bar x(\epsilon)$ the 
$\bar x$ vector implicit in the final iteration of the modified search.  Similarly denote by $V(\epsilon)$ the 
final value of $V$.  Since the
feasible region for (\ref{cdtgen}) is contained in a compact set, we have
that as $\epsilon \rightarrow 0^+$, the sequence of vectors
$\bar x(\epsilon)$  will have an accumulation point which must be an optimal
solution for (\ref{cdtgen}), and the corresponding accumulation point of
the values $V(\epsilon)$ is the optimal value for problem (\ref{cdtgen}).  If on the other hand (\ref{cdtgen}) is infeasible
then for $\epsilon > 0$ small enough the modified binary search will return
infeasibility.  An open question is whether exact (i.e., not weak) feasibility of (\ref{cdtgen}) can be proved, using a weak feasibility oracle with a choice for $\epsilon$ where $\log \epsilon^{-1}$ is polynomially bounded in
the size of the coefficients in (\ref{cdtgen}).

\subsection{Computing explicit solutions in polynomial time}\label{explicit}
Here we will show 
that for each fixed $m$ there is a polynomial-time algorithm that, given 
a system of the form 
\begin{subequations}
\label{basic}
\begin{eqnarray}
&& f_i(x) \ \le \ 0, \quad 1 \le i \le q,  \\ \label{mquad}
&& \norm{x}^2 - 1 \ = \ 0, \label{mqnorm}
\end{eqnarray}
\end{subequations}
where the $f_i(x)$ are quadratic polynomials, and $0 < \epsilon < 1$, either
proves the system is infeasible, or computes an explicit $\epsilon$-feasible solution. Any system of the form (\ref{sys2}) used in the algorithm in Section \ref{algorithm1} 
can be reduced, by scaling, to an equivalent system (\ref{basic})
with $q = O(m)$ and $n$ appropriately redefined.  For $x \in \R^n$ write
$$ \bm{s(x)} \doteq | \norm{x}^2 - 1|.$$
Algorithm C, given next, sequentially computes
values $\hat x_1, \hat x_2, \ldots, \hat x_n $,  such that at termination either the system
(\ref{basic})
is proved infeasible or the vector $\hat x$ is $\epsilon$-feasible.  In preparation for the algorithm we introduce some notation.

\begin{DE}\label{mudef} Given $0 < \delta $ we let $0 < \mu(\delta) < \delta/4$ be any rational
so that whenever $x, y \in \R^n$ are such that $\norm{x} < 2$ and $\norm{x - y} < \mu(\delta)$ then $|f_i(y) - f_i(x)| \le \delta$ for $1 \le i \le q$, 
and $| s(x) - s(y) | \le \delta/2$. Such a value $\mu(\delta)$ can be computed
in polynomial time in the number of bits in the data and $\log \delta^{-1}$.
\end{DE}

\begin{center}
  \textsc{{\bf Algorithm C}}\vspace*{5pt}\\
  \fbox{
    \begin{minipage}{0.9\linewidth}

\noindent {\bf Setup.} Set $\bm{\Delta = \frac{\mu^2(\epsilon/2)}{n}}$. Choose $0 < \rho < \min\{ \mu(\Delta/2), \Delta/n \}$. \\

\noindent {\bf Initialization.} Set $k = 1$.\\

\noindent {\bf Step 1.} Let $z$ denote the vector $(\hat x_1, \ldots, \hat x_{k-1}, z_k, \ldots, z_n)^T$ where for $1 \le j \le k-1$  the $\hat x_j$ are values computed in prior iterations.  Let
$$ Z^{k} \quad \doteq \quad \left\{z \, : \, f_i(z) \, \le \, (k - 1) \Delta \quad \mbox{for $1 \le i \le q$}, \quad s(z) \le (k - 1) \Delta \right\}.$$

\noindent For each index $h = k, k + 1, \ldots, n$ use the algorithm in Section \ref{binsearch} to
produce one of the following two outcomes (1a) and (1b):
\begin{itemize}
\item [{\bf (1a)}] Decide that $Z^k = \emptyset$.  If so, set $k^* = k - 1$, {\bf stop}, and declare
(\ref{basic}) infeasible.
\item [{\bf (1b)}] Decide that $Z^k$ is $\frac{\rho}{2}$-feasible and compute  rationals $P_{k,h}$ and $M_{k,h}$ such that
$$P_{k,h} \ \mbox{is a $\rho/2$-estimate for}  \  \max\{z_h \, : \, z \in Z^k\}, \ \mbox{and} \ M_{k,h} \ \mbox{is a $\rho/2$-estimate for}  \  \min\{z_h \, : \, z \in Z^k\}$$
\end{itemize}

\noindent {\bf Step 2.} Let $\bar h  \doteq \argmax\{ \max\{P_{k,h}, -M_{k,h}\} \, : \, k \le h \le n \}$. By re-indexing, if necessary, assume $k = \bar h$. 
If $P_{k,k} \ge - M_{k,k}$ we obtain  $\hat x_k$ by rounding
$P_{k,k}$ to the nearest integer multiple of $\rho/2$.  Else, 
we obtain  $\hat x_k$ by rounding
$M_{k,k}$ to the nearest integer multiple of $\rho/2$.\\

\noindent {\bf Step 3.} {\bf If} $k < n$ and $\sum_{j = 1}^k \hat x_k^2 \le 1 - \Delta/2 + \rho^2$, set 
$k \leftarrow k+1$ and {\bf go to Step 1}.  {\bf Otherwise}, define $k^* = k$ and set $\hat x_h = 0$ for $k+1 \le h \le n$. and {\bf stop} the procedure.\\
\end{minipage}
  }
\end{center}

\noindent We now analyze Algorithm C.  

\begin{LE} \label{fast} Algorithm C runs in polynomial time. \end{LE}
\noindent {\em Proof.} First we note that 
$\Delta$ and $\rho$ can be computed in polynomial time in the number of bits in the data, and in $\log \epsilon^{-1}$.  Thus
at each iteration $k$ the size of the description of  $Z^k$ is polynomially bounded in the
number of bits in the data and $\log \epsilon^{-1}$. \QED

\begin{LE} \label{central} Let $k \le n$. Suppose that the algorithm 
does not stop at Step 1a of any iteration $k' \le k$, and that $Z^k \neq \emptyset$ .  
There is a vector $\tilde z$, which is $\rho/2$-feasible for $Z^{k}$ with $|\hat x_k - \tilde z_k| \le \rho$.
\end{LE}
\noindent {\em Proof.} Since at iteration $k$ we obtain outcome 1b, 
there is a vector 
$$ \tilde z \ = \ (\hat x_1, \ldots , \hat x_{k-1}, \tilde z_k, \ldots, \tilde z_n)^T,$$
which is $\rho/2$-feasible for $Z^k$, and such that 
$$ | \tilde z_k - \max\{ P_{k,k} , -M_{k,k} \}| \le \rho/2.$$ 
The result follows from the rounding step used to obtain $\hat x_k$ in Step 2. \QED

\begin{LE} \label{iterative}
Suppose system (\ref{basic}) is feasible. Then for any $k \le k^*$,  $Z^{k} \neq \emptyset$.  
\end{LE}
\noindent {\em Proof.}  By induction on $k$.  For $k = 1$ the result follows
since $Z^{1}$ is the set of points feasible for system (\ref{basic}). 
Moreover suppose $Z^{k} \neq \emptyset$ and $k < k^*$. Let $\tilde z$ 
be the vector produced by Lemma \ref{central}, and 
define
$$ \tilde w \ = \ (\hat x_1, \ldots , \hat x_{k-1}, \hat x_{k}, \tilde z_{k+1}, \ldots, \tilde z_n)^T.$$
By Lemma \ref{central}
$\tilde z$ is $\rho/2$-feasible for $Z^{k}$
and $\norm{ \tilde w - \tilde z} \, \le \, \rho \, < \, \mu(\Delta/2)$.
The definition of the $\mu$ function thus implies that $\tilde w$ is $(\rho + \Delta)/2$-feasible for $Z^{k}$. Since $\rho < \Delta$ this implies
$\tilde w \, \in \, Z^{k + 1}$. \QED
\begin{CO} \label{terminfeas} If the algorithm ever stops at Step 1a, system (\ref{basic}) is infeasible.
\end{CO}

\begin{RE} A variation on the proof for Lemma \ref{iterative} shows that if the algorithm does stop at Step 1a, then it does so with $k = 1.$ \end{RE}

\begin{LE} \label{final}  Suppose the algorithm does not stop at Step 1a.  Then at termination the vector $\hat x$ is $\epsilon$-feasible for system (\ref{basic}). \end{LE}
\noindent {\em Proof.}  Let $\tilde z$ be the vector obtained by applying
Lemma \ref{central} at $k = k^*$.  Thus
\begin{subequations} \label{sumdiff}
\begin{eqnarray}
&& \sum_{k = 1}^n \tilde z_k^2 \ \le \ 1 + k^* \frac{\Delta}{2} \ < \ 1 + (n-1) \Delta, \quad \mbox{and} \label{zk*} \\
&& \sum_{k = 1}^{k^*} |\hat x_k - \tilde z_k|^2 \ = \ |\hat x_{k^*} - \tilde z_{k^*}|^2 \le \ \rho^2. \label{sumdiff2}
\end{eqnarray}
\end{subequations}
If $k^* < n$ then since $\sum_{k = 1}^{k^*} \hat x_k^2 \ge 1 -  \Delta/2 + \rho^2$, we have that
(\ref{sumdiff2}) and $\rho \le \mu(\Delta/2)$ imply
\begin{eqnarray}
&& \sum_{k = 1}^{k^*} \tilde z_k^2 \ \ge \ 1 - \Delta/2 - \Delta/2 + \rho^2 = 1 - \Delta + \rho^2, \label{sum2z}
\end{eqnarray}
 and so, by (\ref{zk*}),
$$ \sum_{k = k^* + 1}^n |\hat x_k - \tilde z_k|^2 \ = \ \sum_{k = k^* + 1}^n \tilde z_k^2 \ \le \ n \Delta - \rho^2.$$  
Hence whether $k^* < n$, or not,
\begin{eqnarray}
&& \norm{\hat x - \tilde z}^2 \ \le \ n\Delta \ = \ \mu^2(\epsilon/2) \label{phew}
\end{eqnarray}
(by definition of $\Delta $).  
Now $\tilde z$ is $\rho$-feasible for $V^{k^*}$, and so 
$\tilde z$ is $n \Delta$-feasible 
for (\ref{basic}).  Using $\Delta = \frac{\mu(\epsilon/2)}{n} < \epsilon/(2n)$ (Definition \ref{mudef}) this implies $\tilde z$ is $\epsilon/2$-feasible for (\ref{basic}).  This fact, together with (\ref{phew}) and the definition of the
$\mu$ function now implies that $\hat x$ is $\epsilon$-feasible for (\ref{basic}),
as desired.   \QED\\

\noindent Corollary \ref{terminfeas} and Lemma \ref{final} complete the proof of Theorem \ref{main}.  However we can prove an additional result, 
namely that at termination of the algorithms all entries $\hat x_k$ 
with $k \le k^*$ are nonzero, and ``large.''  

\begin{LE} \label{it2}  Let $k \le k^*$. Then  $|\hat x_k| \ \ge \ [ \frac{\Delta}{4n}]^{1/2} - \rho/2 \ \ge \ \frac{1}{2}[ \frac{\Delta}{4n}]^{1/2}$.
\end{LE}
\noindent {\em Proof.} Since the algorithm has not terminated by iteration
$k -1$, we have $\sum_{j = 1}^{k-1} \hat x_j^2 \le 1 - \Delta/2 + \rho^2$.  Since at iteration
$k$ the algorithm does not stop in Step 1a, there is a 
vector $\breve z$ that is $\rho/2$-feasible for $Z^k$. Thus
\begin{eqnarray}
&& \sum_{j = 1}^{k-1} \hat x_j^2 + \sum_{j = k}^n \breve x_j^2 \ \ge \ 1 - n \rho/2, \label{float}
\end{eqnarray}
As a result, $\sum_{j = k}^n \breve x_j^2$ is ``large,'' i.e.
$ \sum_{j = k}^n \breve x_j^2 \ > \ \Delta/2 - n \rho/2 - \rho^2 > \Delta/4.$
This concludes the proof, since as per Step 2 of the algorithm $|\breve x_k| \ge (\frac{\Delta}{4n})^{1/2}$ . \QED\\

\noindent {\bf Acknowledgment.} This work was partly supported by ONR award N00014-13-1-0042, LANL award ``Grid Science'' and DTRA award  HDTRA1-13-1-0021. \\

\tiny Mon.Feb.23.083927.2015@littleboy\\
\tiny Tue.Nov.26.220720.2013\\
\end{document}